\documentclass[onecolumn]{IEEEtran}

\usepackage{amsmath,amsfonts}
\usepackage{bm}
\usepackage{graphicx}
\begin{document}

\title{Density Estimation on the rotation group using Diffusive wavelets}

\author{Nicolas~Le Bihan,~Julien~Flamant and~Jonathan H.~Manton
\thanks{Nicolas Le Bihan's research was supported by the ERA, European Union, through the International Outgoing Fellowship (IOF GeoSToSip 326176) program of the 7th PCRD}}

\markboth{Journal of Advances in Information Fusion}{Special Issue on Estimation Involving Directional Quantities}

\maketitle




\begin{abstract}
This paper considers the problem of estimating probability density functions on the rotation group $SO(3)$. Two distinct approaches are proposed, one based on characteristic functions and the other on wavelets using the heat kernel. Expressions are derived for their Mean Integrated Squared Errors. The performance of the estimators is studied numerically and compared with the performance of an existing technique using the De La Vall\'ee Poussin kernel estimator. The heat-kernel wavelet approach appears to offer the best convergence, with faster convergence to the optimal bound and guaranteed positivity of the estimated probability density function.

\end{abstract}





\IEEEpeerreviewmaketitle 




\section{Introduction}
Statistics on Lie groups have became more and more popular in the last decade. Applications can be found in shape statistics \cite{fletcher2003gaussian,fletcher2004principal}, medical imaging \cite{fletcher2003statistics}, multiple scattering processes \cite{le2009nonparametric,Said2010}, crystallography \cite{Hielscher2013}, robotics and mechanics \cite{chirikjian2000engineering,Chirikjian2012} and many others. In directional statistics, the abundance of datasets taking values on spheres and homogeneous spaces has also motivated the study of random variables and processes on Lie groups \cite{Mardia}. Amongst all the matrix Lie groups, the most popular one is the rotation group in ${\mathbb R}^3$, {\em i.e} the Special Orthogonal group $SO(3)$. This is due to its predominent use in engineering problems \cite{chirikjian2000engineering}.

Even though the number of engineering challenges including random rotations has grown dramatically in the last years, the problem of probability density estimation for such variables was only considered in details recently in \cite{Hielscher2013}. In parallel, the concept of wavelets on manifolds was transposed to the case of the rotation group, due to its relation to the 2-sphere which has attracted a lot of work since the 90s \cite{schroder1995spherical,antoine1999wavelets,starck2006wavelets}. In \cite{Ebert2011}, authors introduce diffusive wavelets on manifolds. The definition of a wavelet transform on a surface/manifold is conditioned by the possible definition of two operations on the manifold: scaling and translation. While translation is easily defined on Lie groups as it is based on the group action, the definition of scaling is less obvious. In \cite{Ebert2011}, authors have chosen an intrinsic definition for scaling, whereas some extrinsic approaches had been proposed previously in \cite{bogdanova2005stereographic}. The difference in these definitions resides in the way the mother wavelet is scaled: it is either firstly projected in the tangent plane before scaling and back projection (extrinsic); or scaled on the manifold directly (intrinsic). In this paper, we will make use of the intrinsic approach and study the ability of diffusive wavelets to define interesting estimators for densities on $SO(3)$. 

The specificity of diffusive wavelets is that they are based on the heat kernel. As this kernel can be defined on manifolds, the diffusive wavelet approach overcomes the problem of ``scaling'' on manifold. Note that this definition problem was already pointed out by many authors \cite{bogdanova2005stereographic,	antoine1999wavelets} when defining wavelets on the 2-sphere for example. 

In this paper, we present an estimation technique for densities on $SO(3)$ based on the diffusive wavelet transform. Using wavelet estimators (linear or thresholded) is well known in non parametric estimation \cite{Tsybakov2009}. We propose the use of the linear\footnote{Linear wavelet estimation refers to the standard estimation through wavelet coefficients estimation, as opposed to thresholded wavelet coefficients estimator sometimes refered as nonlinear \cite{Tsybakov2009}} wavelet estimator to the case of $SO(3)$-valued random variables and give some of its properties. In particular, the wavelet-based estimator is a characteristic kernel estimator. A comparison with other types of estimators (characteristic function and kernel) is provided. 

The rest of the paper is organized as follows: Section \ref{sec:CharaF_SO3} is dedicated to the presentation of harmonic analysis on $SO(3)$ and various concepts about random variables on this Lie group. In Section \ref{sec:DiffWavSO3}, the diffuse wavelet formalism is introduced for the case of the rotation group. It is then used in Section \ref{sec:estim} to propose a probability density function estimator. Finally, Section \ref{sec:simu} presents simulation results and comparison of the wavelet estimator with two other estimators.



	
\section{Characteristic functions for random variables on $SO(3)$}
\label{sec:CharaF_SO3}
We first review some basic concepts on the rotation group $SO(3)$ and random variables taking values on this well known Lie group. The concept of characteristic function for $SO(3)$-valued random variables is of importance as it can be identified as the ``Fourier transform'' of the probability density of the random variable. It will also be of use in Section \ref{sec:estim} to provide a probability density estimator. The results presented here can be found in various textbooks, for example \cite{Altmann2005,Dym1972}.

%
%

\subsection{The rotation group $SO(3)$}
The set of rotations in the 3D space forms a compact Lie group denoted $SO(3)$. An element $x \in SO(3)$ can be parametrized in several ways \cite{Altmann2005}. Using the so-called ZYZ convention of Euler angles parametrization, any element $x \in SO(3)$ can be associated to a matrix $R^x=R^x(\varphi,\theta,\psi)$ with $0 \leq \varphi,\psi < 2\pi$ and $0 \leq \theta \leq \pi$. With this convention, $R^x(\varphi,\theta,\psi)$ takes the form: 
\begin{align*}
R^x(\varphi,\theta,\psi) = &
\left[
\begin{array}{ccc}
\cos \varphi & -\sin \varphi & 0 \\
\sin \varphi & \cos\varphi & 0 \\
0 & 0 & 1
\end{array}
\right]
\left[
\begin{array}{ccc}
\cos \theta & 0 & \sin \theta \\
0 & 1 & 0 \\
-\sin \theta & 0 & \cos\theta
\end{array}
\right]
\\
& \times 
\left[
\begin{array}{ccc}
\cos \psi & -\sin\psi & 0 \\
\sin \psi & \cos\psi & 0 \\
0 & 0 & 1
\end{array}
\right]
\end{align*}

It is also possible to parametrize an element of $SO(3)$ in terms of its \emph{rotation axis} and its \emph{rotation angle} \cite[Chap. 3]{Altmann2005}. The rotation axis $\eta$ is a unit vector in ${\mathbb R}^3$, \emph{i.e.} $\eta \in {\cal S}^2$. The rotation angle, denoted $\omega(x)$, is given by:
\begin{equation}
\cos \omega(x)= \frac{Tr(R^x) - 1}{2}
\label{eq:rotangleS03}
\end{equation} 
Note that in this parametrization, the angle takes values in: $-\pi < \omega(x) \leq \pi$. It is also possible to express this angle $\omega(x)$ in terms of the Euler angles:
\begin{equation}
\omega(x)= 2 \arccos \left( \cos \frac{\theta}{2} \cos \frac{\varphi+\psi}{2}\right)
\end{equation} 
The \emph{rotation angle} is of particular interest as it is a metric on $SO(3)$. In fact, one can define the distance between $x \in SO(3)$ and $y \in SO(3)$, denoted $d(x,y)$, as:
\begin{equation}
d(x,y)=|\omega(yx^{-1})| = \arccos \left(\frac{Tr(R^y(R^{x})^{-1}) - 1)}{2}\right)
\label{eq:geodistSO3}
\end{equation}   
In the sequel, we will make use of the notation ${\bm R(\omega(x),\eta)}$ for an element $x \in SO(3)$, keeping in mind that:
\begin{equation}
{\bm R(\omega(x),\eta)} = \exp(\omega(x){\bm M}) \quad \text{with} \quad {\bm M}=\left(\begin{matrix} 0 & -\eta_3 & \eta_2 \\ \eta_3 & 0 & -\eta_1 \\ -\eta_2 & \eta_1 & 0 	\end{matrix}\right)
\end{equation}
where $\eta_1,\eta_2,\eta_3$ are the components of the vector $\eta$, {\em i.e.} $\eta=\left[\eta_1,\eta_2,\eta_3\right]^T$ and $\exp(.)$ the matrix exponential \cite{Altmann2005}. 

\subsection{Fourier series on $SO(3)$}
We now consider the Fourier series expansion of functions taking values on the rotation group. Consider the set of square integrable functions $L^2(SO(3),{\mathbb R})$. By the Peter-Weyl theorem \cite{Dym1972}, a function $f \in L^2(SO(3),{\mathbb R})$ can be expressed as:
\begin{equation}
f(x)= \sum_{\ell\geq 0}\sum_{n = -\ell}^{+\ell}\sum_{m = -\ell}^{+\ell} (2\ell+1) \hat{f}_{nm}^{\ell} D_{nm}^{\ell}(x)
\label{eq:FourierSeriesSO3}
\end{equation} 
for $x \in SO(3)$ and where $D_{nm}^{\ell}(x)$ are the Wigner-D functions \cite{Barut1986} evaluated at position $x$. This infinite series expansion over $\ell$ has matrix coefficients $\hat{f}^{\ell}$ of dimension $(2\ell+1)\times(2\ell+1)$ with elements $\hat{f}_{nm}^{\ell}$. The matrix entries $\hat{f}_{nm}^{\ell}$ are obtained by projection of $f$ on the Wigner-D functions:
\begin{equation}
\hat{f}_{nm}^{\ell} =\left<f,D_{nm}^{\ell}\right>_{\text{\tiny{$SO(3)$}}} = \int_{SO(3)} f(x) \overline{D_{nm}^{\ell}}(x)  d\mu 
\end{equation} 
where $\left<f,h \right>_{\text{\tiny{$SO(3)$}}}$ is the scalar product on $L^2(SO(3),{\mathbb R})$ and $d\mu$ the bi-invariant Haar measure on $SO(3)$, {\em i.e.} $d\mu=(8\pi^2)^{-1}\sin\theta d\varphi d\theta d\psi$ when using the ZYZ Euler angle parametrization for elements $x \in SO(3)$. 	
We also mention that the corresponding norm is:
\begin{equation}
\left\|f\right\|_2=\sqrt{\left<f,f \right>_{\text{\tiny{$SO(3)$}}}}
\end{equation}
With the chosen parametrization of $SO(3)$, the previously introduced Wigner-D functions $D_{nm}^{\ell}$ take the form:
\begin{equation}
D_{nm}^{\ell}\left(\varphi,\theta,\psi\right)= e^{-{\tt i}n\varphi} P_{nm}^{\ell}(\cos\theta) e^{-{\tt i}m\psi}
\end{equation}
where $P_{nm}^{\ell}(\cos\theta)$ are the generalized Legendre polynomials. We have used the fact that the Wigner-D functions $D_{nm}^{\ell}$ form a complete set of orthogonal functions:
\begin{equation}
\begin{array}{rcl}
\displaystyle{\left<D_{nm}^{\ell},D_{n'm'}^{\ell'}\right>_{\text{\tiny{$SO(3)$}}}} 
 & =  & \displaystyle{\frac{1}{\left(2\ell + 1\right)} \delta_{nn'}\delta_{mm'}\delta_{\ell\ell'}} 
\end{array}
\end{equation}
This shows clearly that the set of functions $\left\{\sqrt{2\ell+1}D^{\ell}_{nm}, \ell \geq 0, -\ell\leq n,m \leq \ell \right\}$ form an orthonormal basis for functions in $L^2(SO(3),{\mathbb R})$, allowing the decomposition of functions taking values on the rotation group using this basis. In the sequel, we may refer to $\hat{f}_{nm}^{\ell}$ (and abusively to $\hat{f}^{\ell}$) as the {\em Fourier coefficients} of $f$. 

%
%
\subsection{Parseval identity}
It is also well known that the Parseval identity holds for functions $f \in L^2(SO(3),{\mathbb R})$. Using our normalization convention, the following is true:
\begin{equation}
\begin{array}{rcl}
\left\|f\right\|_2^2 & = & \displaystyle{\left<f,f\right>_{\text{\tiny{$SO(3)$}}}} \\
& = & \displaystyle{\sum_{\ell \geq 0} \sum_{n=-\ell}^{\ell} \sum_{m=-\ell}^{\ell} (2\ell +1) \left|\hat{f}^{\ell}_{nm}\right|^2} 
\end{array}
\label{eq:parsevalSO3}
\end{equation} 
The left hand side of equation (\ref{eq:parsevalSO3}) is commonly denoted as the {\em energy} and the Parseval identity simply states that the energy of $f$ consists in the infinite sum of its modulus squared {\em Fourier coefficients}.  

\subsection{Zonal functions on $SO(3)$}
\label{sec:zonalfunctSO3}
Zonal functions, sometimes also called \emph{conjugate invariant} functions, are radially symmetric functions with center $g_0=e \in SO(3)$ where $e$ is the identity element in $SO(3)$. See for example \cite{Hielscher2007} for more details on radially symmetric fonctions on $SO(3)$. They will be at the heart of the diffusive wavelets construction in Section \ref{sec:DiffWavSO3}. A function $f:SO(3) \rightarrow {\mathbb R}$ is called \emph{zonal} iff it satisfies $\forall x,y \in SO(3)$:
$$ f(yxy^{-1})=f(x)$$
Equivalentlty, $f:SO(3)\rightarrow {\mathbb R}$ is \emph{zonal} iff $f(x)=f(x')$ $\forall x,x' \in SO(3)$ such that $\omega(x)=\omega(x')$, {\em i.e.} $x$ and $x'$ have the same {\em angle}. In fact, it simply means that $f$, evaluated at $x \in SO(3)$, only depends on the rotation angle of $x$ introduced in (\ref{eq:rotangleS03}). It is known \cite{Hielscher2007,Hielscher2013} that the subspace of zonal functions is spanned by functions $\chi^{\ell}$, $\ell \in {\mathbb N}$, called the {\em characters} of $SO(3)$, and given as:
\begin{eqnarray}
\chi^{\ell}(x) & = & \sum_{n=-\ell}^{\ell} D^{\ell}_{nn}(x) = \frac{\sin \left(\left( \ell + \frac{1}{2}\right)\omega(x)\right)}{\sin \frac{\omega(x)}{2}}\\
 & = & {\mathcal U}_{2\ell}\left(\cos\frac{\omega(x)}{2}\right)
\end{eqnarray}
where ${\mathcal U}_{2\ell}(.)$ are the Chebychev polynomials of second kind and of (even) degree $2\ell$. For these (even degree) polynomials, the following orthogonality relation stands:
\begin{equation}
 \int_{-1}^1 {\mathcal U}_{2\ell}(t) {\mathcal U}_{2\ell'}(t) \sqrt{1-t^2}dt = \frac{\pi}{2}\delta_{\ell\ell'}
\end{equation} 
In terms of the characters $\chi^{\ell}$, and using the notation $\omega(x)=\omega$ for simplicity, this integral becomes:
\begin{equation}
\int_0^{2\pi}\chi^{\ell}(\omega)\chi^{\ell'}(\omega)\sin^2\left(\frac{\omega}{2}\right)d\omega = \pi \delta_{\ell \ell'}
\end{equation}
where the link with the orthogonality relation for ${\mathcal U}_{2\ell}$ comes with $t=\cos\frac{\omega}{2}$. As a consequence, the characters $\chi^{\ell}$ fulfill the following orthogonality relation:
\begin{equation}
\left< \chi^{\ell}, \chi^{\ell'} \right>_{\text{\tiny{$SO(3)$}}} =   \frac{1}{4\pi^2} \int_{{\mathcal S}^2}d\eta \int_0^{2\pi} \chi^{\ell}(\omega) \chi^{\ell'}(\omega) \sin^2\left(\frac{\omega}{2}\right)d\omega =  \delta_{\ell\ell'}
\end{equation}
where we used the expression of the Haar measure in terms of axis $\eta$ and rotation angle $\omega$ for $SO(3)$. With this convention the Haar measure takes the form $d\xi=\frac{1}{4\pi^2}d\eta \sin^2(\frac{\omega}{2})d\omega$. As a consequence, a zonal function $f \in L^2(SO(3),{\mathbb R})$ has a Fourier series expansion that can be written like:
\begin{equation}
f(x)=\sum_{\ell \geq 0} (2 \ell +1) \hat{F}^{\ell} \chi^{\ell}(x)
\label{eq:FourierZonal}
\end{equation}
where its Fourier coefficients $\hat{F}^{\ell}$ are denoted in capital to emphasize the fact that $f$ is zonal. These Fourier coefficients $\hat{F}^{\ell}$ are simply given by:
\begin{equation}
\hat{F}^{\ell} =  \frac{1}{(2 \ell +1)} \left< f,\chi^{\ell}\right>_{\text{\tiny{$SO(3)$}}} = \frac{1}{\pi} \frac{1}{(2 \ell +1)}\int_0^{2\pi} f(\omega) \chi^{\ell}(\omega) \sin^2 \left(\frac{\omega}{2}\right) d\omega
\end{equation}
The Parseval identity for a zonal function $f$ now reads:
\begin{equation}
\|f\|^2_2= \sum_{\ell \geq 0} \left(2\ell+1\right)^2\left|\hat{F}^{\ell}\right|^2
\label{eq:ParsevalZonal}
\end{equation}
Zonal functions will be used in Section \ref{sec:DiffWavSO3} in the definition of wavelets on $SO(3)$.

\subsection{Convolution}
A very important feature of Fourier transformation is its behaviour with respect to the convolution product. In the case of functions taking values on $SO(3)$, the "famous" relation still holds. First, given $f,h \in L^1(SO(3),{\mathbb R})$, their convolution product is defined as: 
\begin{equation}
\left(f*h\right)(x) = \int_{\text{\tiny{$SO(3)$}}} f(g)h(g^{-1}x) d\mu
\end{equation}
where the group operation stands naturally for translation. Now, if the Fourier coefficients of $f$ and $h$ are respectively $\hat{f}^{\ell}$ and $\hat{h}^{\ell}$ (in matrix format), then the following stands:
\begin{equation}
\widehat{f*h}^{\ell}=\hat{f}^{\ell}\hat{h}^{\ell}
\end{equation}
Note that the right-hand side of the equation is a matrix product.

\subsection{Characteristic function of $SO(3)$-valued random variables}
The characteristic function of a random variable is the Fourier transform of its probability transform. This well known result extends to random variables on the rotation group, thanks to the results on Fourier series expansion introduced in the previous Section. We now give some definitions and properties for characteristic function of random variables on $SO(3)$ as it will be used to define a probability density function estimator in Section \ref{sec:estim}. 

Consider the case of a random variable $X \in SO(3)$ with density $f$. Such random variables $X$ taking values on the rotation group $SO(3)$ can simply be thought of as random rotation matrices of dimension $3 \times 3$, classically parametrized by Euler angles. Alternatively, one can think of these random variables as unit quaternions from the upper hemisphere of ${\mathcal S}^3$. 

\subsubsection{Definition}
Given a rotation random variable $X$ with density $f_X$, the sequence $\Phi_X=\left\{\Phi_X(\ell)\right\}_{\ell \geq 0}$ of $(2\ell + 1) \times (2\ell +1)$ matrices given by:
\begin{equation}
\Phi_X(\ell)={\mathbb E}\left[\overline{D^{\ell}}(X)\right] 
\end{equation}
is the \emph{characteristic function} of $X$. The elements of the matrix $\Phi_X(\ell)$ are denoted $\Phi_{nm}^{\ell}$ and they read as:
\begin{equation}
\Phi_{nm}^{\ell}=\int_{SO(3)}f_X(x)\overline{D^{\ell}_{nm}}(x) d\mu
\end{equation}
Thus the density $f$ has the following Fourier series expansion:
\begin{equation}
f_X(x)= \sum_{\ell \geq 0}\sum_{n,m=-\ell}^{\ell} (2\ell + 1)\Phi_{nm}^{\ell} D^{\ell}_{nm}(x)
\end{equation}

One can see by comparison with Equation (\ref{eq:FourierSeriesSO3}) that the characteristic function of $f_X$ consists is the Fourier transform of its density $f_X$, \emph{i.e.} $\Phi_X^{\ell}=\hat{f}_X^{\ell}$.

\subsubsection{Basic properties}
The following properties can easily be verified:
\begin{itemize}
\item Given two rotation random variables $X$ and $Y$, then:
$$X=Y \quad \text{iff} \Phi_X=\Phi_Y$$
\item If $X$ and $Y$ are two independent rotation random variables and $Z=XY$, then:
$$\Phi_Z^{\ell}=\Phi_X^{\ell}\Phi_Y^{\ell}$$
\item A rotation random variable $U \in SO(3)$ is uniformly distributed iff:
$$\Phi_U^{\ell} = 0 \quad \forall \quad \ell > 0$$
\item Consider $n$ \emph{i.i.d.} rotation random variables $X_1,X_2,\ldots,X_n$, then the random variable consisting in the accumulated products of the $X_n$, denoted $Y=X_1X_2 \ldots X_n$ has the following characteristic function: 
$$\Phi_Y^{\ell}=\left[\Phi_{X}^{\ell}\right]^n$$ 
where $\Phi_{X}^{\ell}$ is the characteristic function shared by the $X_n$.
\end{itemize}

\subsubsection{Zonal invariance}
	A rotation random variable $X$ with density $f_X$ is called \emph{zonal invariant} if:
	$$X \stackrel{d}{=} R X R^{-1} \quad \text{for all} \quad R \in SO(3)$$
	Now, if $X$ is \emph{zonal invariant}, then its characteristic function is: 
	$$\Phi_X^{\ell} = a_{\ell} I_{\ell}$$ 
	where $a_{\ell} \in {\mathbb R}$ and $I_{\ell}$ is the $(2\ell+1) \times (2\ell+1)$ identity matrix. As a consequence, if $X$ is \emph{zonal invariant}, then its density $f_X$ takes the form:
	\begin{equation}
	f_X(x) = \sum_{\ell \geq 0}  \, (2\ell+1) a_{\ell} \, \chi^{\ell}(x)
	\end{equation} 
	as detailed previously when considering zonal functions on $SO(3)$ in \ref{sec:zonalfunctSO3}.

This last expression shows that the characteristic function of {\em zonal invariant} random rotation variables are scalar coefficients $a_{\ell}$, as opposed to matrix coefficients for random variables with no symmetries. 

%



\section{Diffusive wavelets on $SO(3)$}
\label{sec:DiffWavSO3}
We now introduce diffusive wavelets on the rotation group. They will be used in Section \ref{sec:estim} to propose a probability density estimator. Recently, Ebert and Wirth \cite{EbertThesis2011} introduced \emph{diffusive wavelets} on groups. We detail in the sequel the special case of the rotation group in 3D, {\em i.e.} $SO(3)$. 

\subsection{Heat wavelet family on $SO(3)$}
\label{subsec:HeatWaveletSO3}
Here, we follow the construction given by Ebert \cite{EbertThesis2011}. First recall that on $SO(3)$ the heat kernel is given by:
\begin{equation}
\kappa_{\rho}(x)=\sum_{\ell \geq 0}(2\ell + 1 ) e^{- \ell (\ell +1) \rho} \chi^{\ell}(x)
\end{equation}
for $x \in SO(3)$, and where $\chi^{\ell}$ are the irreducible characters of $SO(3)$ introduced previously in Section \ref{sec:CharaF_SO3}. Note that the Fourier series expansion of $\kappa_{\rho}$ exhibits the fact that $\kappa_{\rho}$ is a zonal function with Fourier coefficients $\hat{\kappa}_{\rho}^{\ell}=e^{- \ell (\ell +1) \rho}$. 

Now, we introduce some of the properties of the heat kernel. First, $\kappa_{\rho}(x)$ is a \emph{diffusive approximate identity} as it fullfils the following properties: 

\begin{itemize}
\item $\left\|\hat{\kappa}_{\rho}^{\ell}\right\| \leq C \quad \forall \ell \quad \rho \in {\mathbb R}^+$
\item  $\lim_{\rho \rightarrow 0} \hat{\kappa}^{\ell}_{\rho}= Id \quad \forall \ell$
\item $\lim_{\rho \rightarrow \infty} \hat{\kappa}^{\ell}_{\rho}= 0 \quad \forall \ell > 0$
\item $-\frac{\partial}{\partial \rho} \hat{\kappa}^{\ell}_{\rho}$ is a positive matrix for all $\rho > 0$ and $\ell \geq 0$
\end{itemize} 
where $Id$ denotes the identity operator. In fact, $\left\{\kappa_{\rho},\rho > 0\right\}$ defines a family of kernels associated to an approximate identity (\emph{i.e.} an admissible semi-group\footnote{A continuous family $\left\{D_\rho,\rho>0\right\}$ of operators on $L^2(SO(3))$ is called an admissible semi-group if:
\begin{itemize}
\item $D_\rho$ is a bounded operator independent of $\rho$
\item $\lim_{\rho \rightarrow 0} = Id$
\item $D_\rho$ is postive $\forall \rho$
\item $D_{\rho_1}D_{\rho_2} = D_{\rho_1+\rho_2}$
\end{itemize}
Now, an admissible semi-group that can be expressed as a convolution operator, \emph{i.e.} for the semi-group $\left\{D_\rho,\rho>0\right\}$ there exits a family of kernels $\left\{\kappa_\rho, \rho>0\right\}$ such that $D_\rho(f)=f*\kappa_\rho$. Such operator is also called \emph{approximate identity with kernel $\kappa_\rho$}.} with kernel $\kappa_{\rho}$). It means that: 
$$
\kappa_{\rho} * h \underset{\rho \longrightarrow 0}{\longrightarrow} h \quad \forall h \in L^2(SO(3)) 
$$  

As detailed in \cite{EbertThesis2011}, a family of wavelets corresponding to the heat kernel on $SO(3)$ is consequently of the form:
\begin{equation}\label{eq:wavheatso3}
\displaystyle{\Psi_\rho(x)=\frac{1}{\sqrt{\alpha(\rho)}} \sum_{\ell \geq 0} (2 \ell +1) \sqrt{\ell (\ell +1)} e^{- \frac{\ell (\ell +1)}{2} \rho} \chi^{\ell}(x)  
}
\end{equation}
with $\alpha(\rho)$ a normalizing factor to be detailed below. Rephrasing this equation in terms of Fourier transform, noticing that the heat wavelet family is zonal, one can write that:
\begin{equation}
\Psi_\rho(x) = \sum_{\ell \geq 0} (2\ell +1) \hat{\Psi}_\rho \chi^{\ell}(x)
\end{equation}
where the Fourier coefficients of $\Psi_\rho$ are: 
\begin{equation}
\hat{\Psi}_\rho = \frac{1}{\sqrt{\alpha(\rho)}} \sqrt{\ell (\ell +1)} e^{- \frac{\ell (\ell +1)}{2} \rho}
\end{equation}
As can be seen in Equation (\ref{eq:wavheatso3}), the choice for the normalization coefficients $\alpha(\rho)$ has to be made. In order to have a normalized/unitary wavelet family, we can use the Parseval identity to choose $\alpha(\rho)$. In fact, by (\ref{eq:ParsevalZonal}), we have that: 
\begin{equation}
\|\Psi_\rho\|^2_2 = \sum_{\ell \geq 0} (2\ell + 1)^2 \left|\hat{\Psi}_\rho\right|^2
\end{equation}
and imposing that $\|\Psi_\rho\|^2_2=1$ $\forall \rho$, one obtains that:
\begin{equation}
\displaystyle{\alpha(\rho)= \sum_{\ell \geq 0} (2\ell + 1)^2 \ell(\ell+1) e^{-\ell(\ell+1)\rho}}
\end{equation}
This choice for $\alpha(\rho)$ will be made throughout the rest of the paper. Notice also that the heat wavelet family  $\Psi_\rho \in L^2(SO(3))$ is a \emph{diffusive wavelet family} as it satisfies the \emph{admissibility condition}:
\begin{equation}
\kappa_\rho(x)  =  \displaystyle{\int_{\rho}^{+\infty} \left(\check{\Psi}_{t} * \Psi_{t}\right)(x) \alpha(t)dt}  =  \displaystyle{\int_{\rho}^{+\infty} \int_{SO(3)} \check{\Psi}_{t}(g) \Psi_{t}(g^{-1}x) d\mu(g) \alpha(t)dt}
\label{eq:heatkernelSO3}
\end{equation}
in which $\kappa_{\rho}(x)$ is an diffusive approximate convolution identity and where we used the notation $\check{\Psi}_{\rho}(g)=\overline{\Psi_\rho(g^{-1})}$. The fact that $\kappa_{\rho}(x)$ is an approximate convolution identity ensures that the wavelet transform can be inverted. Note that the approximate convolution identity can be expressed using the Fourier coefficients. In the case of $\kappa_\rho$, it means that:
$$\lim_{\rho \rightarrow 0} \hat{\kappa}_{\rho}^{\ell}=1 \quad \forall \ell $$
This means that an approximate convolution identity is characterized by constant Fourier coefficients in the limit. 

With the diffusive wavelets introduced, we can now introduce the wavelet transform for functions on the rotation group.

\subsection{Wavelet transform on $SO(3)$}
\label{subsec:WTSO3}
Recall that we are interested in estimating probability density functions of $SO(3)$-valued random variables using the diffusive wavelet transform. A probability density function $f$ belongs to $L^2(SO(3),{\mathbb R}) \cap L^1(SO(3),{\mathbb R})$ with the aditionnal condition that $\displaystyle{\int_{SO(3)} f d\mu}=1$ and that it is non-negative. The diffusive wavelet transform for such densities is as follows. Given a diffusive wavelet family $\Psi_{\rho} \in L^1(SO(3))$ as defined in \ref{subsec:HeatWaveletSO3}, then the Wavelet Transform (WT) of $f \in L^2(SO(3),{\mathbb R})\cap L^1(SO(3),{\mathbb R})$ is:
$$
\begin{array}{rccc}
WT: & f & \longrightarrow & WT_f \\
 & SO(3) & & {\mathbb R}^+ \times SO(3) 
\end{array}
$$
with the following expression:
\begin{equation}
WT_f(\rho,g)= \left( f * \check{\Psi}_{\rho} \right) (g) = \int_{SO(3)} f(x)\check{\Psi}_{\rho}(x^{-1}g) d\mu(x)
\end{equation}
and where we made use of the notation $\check{\Psi}(x)=\overline{\Psi(x^{-1})}$. Using scalar product on the rotation group, this expression can be written like:
\begin{equation}
WT_f(\rho,g)= \int_{SO(3)} f(x) \overline{\Psi_{\rho}(g^{-1}x)} d\mu(x) = \left<f, T^*_g \Psi_\rho \right>_{\text{\tiny{$SO(3)$}}}
\label{eq:wavelso3}
\end{equation}
with $T^*_g$ the following operator: $T^*_g:\Psi \rightarrow \Psi(g^{-1}.)$. Equation (\ref{eq:wavelso3}) is a very general definition of wavelet transform on $SO(3)$. In the sequel, we will only make use of the heat wavelet family given in (\ref{eq:wavheatso3}) to define the heat wavelet transform. One of the interesting properties of the wavelet transform is that it is invertible. The density $f$ can thus be reconstructed in the following way:
\begin{equation}
\begin{array}{rcl}
f(x) & = & \displaystyle{\int_{{\mathbb R}^+}\int_{SO(3)} WT_f(t,g) \Psi_{t}(g^{-1}x) d\mu(g) \alpha(t)dt} \\
 & = & \displaystyle{\int_{{\mathbb R}^+}\int_{SO(3)} \left( f * \check{\Psi}_{t}\right)(g) \Psi_{t}(g^{-1}x) d\mu(g) \alpha(t)dt}\\
 & = & \displaystyle{\int_{{\mathbb R}^+} \left(f * \check{\Psi}_{t} * \Psi_{t} \right)(x) \alpha(t)dt} \\
 & = & \displaystyle{ f * \int_{{\mathbb R}^+} \left( \check{\Psi}_{t} * \Psi_{t} \right)(x) \alpha(t)dt} \\
 & = & \displaystyle{ f * \int_{\rho \rightarrow 0}^{+\infty} \left( \check{\Psi}_{t} * \Psi_{t} \right)(x) \alpha(t)dt} \\
 & = & \displaystyle{ \lim_{\rho \rightarrow 0} \left(f*\kappa_{\rho}\right)(x) }\\
 & = & \displaystyle{f(x)} 
\end{array}
\label{eq:invWTSO3}
\end{equation}
where the last equality is obtained thanks to the fact that $\kappa_t(x)$ is a convolution identity. This can be verified in the Fourier domain, when denoting $\hat{F}^{\ell}$ the Fourier coefficients of $f$, by:
\begin{equation}
\begin{array}{rcl}
 \displaystyle{\lim_{\rho \rightarrow 0} \left( f*\kappa_{\rho} \right)(x)} & = & \displaystyle{\lim_{\rho \rightarrow 0} \sum_{\ell \geq 0} (2\ell +1) e^{-\ell(\ell+1)\rho} \hat{F}^{\ell} \chi^{\ell}(x)}\\
  & = & \displaystyle{\sum_{\ell \geq 0} (2\ell +1) \hat{F}^{\ell} \chi^{\ell}(x)}\\
  & = & \displaystyle{f(x)}
\end{array}
\end{equation}

An interesting property of the heat wavelet transform as defined above is that it is {\em unitary}. This property reads 
\begin{equation}
\left<WT_{f_1}(\rho,g),WT_{f_2}(\rho,g)\right>_{\text{\tiny{${\mathbb R}^{+} \times SO(3)$}}} = \left< f_1,f_2\right>_{\text{\tiny{$SO(3)$}}}
\end{equation}
where we use the notation $\left<.,.\right>_{\text{\tiny{${\mathbb R}^{+} \times SO(3)$}}}$ for the scalar product between wavelet transforms. This unitary condition makes possible to compare density in the wavelet domain for example.

We have introduced the Heat wavelet transform (also called diffusive wavelet transform) for probability densities of rotation random variables. In the next Section, we will make use of this transform to define an estimator for densities on $SO(3)$. 




\section{Estimation}
\label{sec:estim}
In this section, we consider the problem of estimating the probability density function of a $SO(3)$-valued random variable given an independent sample of size $K$: $\left\{X_1,X_2,\ldots,X_K\right\}$. After presenting kernel estimators as defined in \cite{Hielscher2013}, we show that the characteristic function estimator and the heat wavelet estimator are actually kernel estimators. We provide the MISE of each of them and explain their differences. 

%
%

\subsection{Kernel estimators on $SO(3)$}
First, we recall some of the results given in \cite{Hielscher2013} for kernel estimators of probability density functions on $SO(3)$. 

\subsubsection{Definition}
A Kernel estimator, with Kernel $\Xi \in L^2(SO(3))$ has the following expression:
\begin{equation}
\zeta_K (x) = \frac{1}{K}\sum_{k=1}^K \Xi(X_k^{-1}x)
\end{equation}
This definition is the extension of the classical kernel estimator known for densities of random variables on the line (see \cite{Tsybakov2009} for details). Once again, we emphasize that translation is made through the group action. Before introducing different types of kernels $\Xi(.)$, we introduce how to study their estimation performances.

\subsubsection{Bias and variance}
\label{subsec:BVK}
In order to characterize the Kernel estimator, we provide here the expression of its Mean Integrated Square Error (MISE). This expression is general for kernel estimators and will be of use later to analyze the behaviour of the characteristic function and the wavelet estimators.

In the sequel, we will make use of the following notation in order to distinguish the different density estimator. Every estimator, based on a simple sample of size $K$, will be denoted $\zeta^{\Delta}_K(.)$ with the superscript $\Delta$ taking the following values: $\Delta=\alpha$ for the ``general'' kernel estimator, $\Delta=\beta$ for the characteristic function estimator and $\Delta=\gamma$ for the Heat wavelet estimator. Associate quantities will exhibit the $\alpha,\beta,\gamma$ values when needed.

As known in the classical case \cite{Tsybakov2009} and as given for the $SO(3)$ case in \cite{Hielscher2013}, the MISE of the Kernel estimator $\zeta^{\alpha}_K (x)$ is made of a bias and a variance term: 
\begin{equation}
MISE\left(\zeta^{\alpha}_K (x)\right) = \left\|f - {\mathbb E}\left[\zeta^{\alpha}_K\right]\right\|_2^2 + {\mathbb E}\left[ \left\| \zeta^{\alpha}_K (x) - {\mathbb E}\left[\zeta^{\alpha}_K (x)\right] \right\|_2^2 \right]
\end{equation}
with the following fact:
\begin{equation}
\displaystyle{ {\mathbb E}\left[\zeta^{\alpha}_K (x)\right]} = \displaystyle{{\mathbb E}\left[ \frac{1}{K}\sum_{k=1}^K \Xi(X_k^{-1}x) \right]} =  \displaystyle{ \frac{1}{K}\sum_{k=1}^K {\mathbb E}\left[\Xi(X_k^{-1}x) \right] } = \displaystyle{ \frac{1}{K}\sum_{k=1}^K \int_{SO(3)} \Xi(y^{-1}x)f(y)d\mu(y) } = \displaystyle{ \left(\Xi * f \right)(x)}
\end{equation}
This means that the mean value of the Kernel consists in the convolution of the density with the kernel. One can also express the bias $b_{\alpha}$ and the variance $\sigma^2_{\alpha}$ as follows: 
\begin{equation}
b^2_{\alpha}=\left\|f - {\mathbb E}\left[\zeta^{\alpha}_K\right]\right\|_2^2 = \left\|f - f * \Xi \right\|_2^2
\end{equation}
and
\begin{equation}
\sigma^2_{\alpha}={\mathbb E}\left[ \left\| \zeta^{\alpha}_K (x) - {\mathbb E}\left[\zeta^{\alpha}_K (x)\right] \right\|_2^2 \right] = \frac{1}{K} \left( \left\| \Xi \right\|_2^2 - \left\| \Xi * f \right\|_2^2\right)
\end{equation}
Also, it is interesting to note that the following property holds:
\begin{equation}
\lim_{K \to \infty} \zeta^{\alpha}_K (x) = \left(\Xi * f\right)(x)
\end{equation}
As a consequence, the MISE for the kernel estimator is:
\begin{equation}
MISE(\zeta^{\alpha}_K) = \left\|f - f * \Xi \right\|_2^2 + \frac{1}{K} \left( \left\| \Xi \right\|_2^2 - \left\| \Xi * f \right\|_2^2\right)
\end{equation}

This general expression is useful for the study of kernel estimator. Using results from representation theory (Fourier series expansion), it is also possible to write the MISE in terms of the Fourier coefficients of $\Xi$ and $f$. Remembering the Parseval identity and the convolution property of the Fourier series expansion on $SO(3)$, one gets, just like in \cite{Hielscher2013} but with a slight $(2\ell+1)$ factor due to our normalization choice, the following representation:
\begin{equation}
\displaystyle{MISE(\zeta^{\alpha}_K) = \sum_{\ell \geq 1} \left( (2\ell +1) \left.\hat{f}^{l}\right.^2  (1 - \hat{\Xi}^{\ell})^2  + \frac{(2\ell +1)}{K} \left.\hat{\Xi}^{\ell}\right.^2 \left(1 - \left.\hat{f}^{l}\right.^2\right)\right)}
\label{eq:MISE_K}
\end{equation}
where the following notation was used for the Fourier coefficients of the probability density function $f$:
\begin{equation}
\left.\hat{f}^{l}\right.^2 = (2\ell +1)^{-1} \sum_{n,m=-\ell}^{\ell}\left|\hat{f}^{\ell}_{n,m}\right|^2
\end{equation}
This expression of the MISE of a kernel estimator on $SO(3)$ will be used with a specific kernel, namely the De La Vall\'ee Poussin kernel, in Section \ref{sec:simu} for comparison purposes with the characteristic function estimator and the wavelet estimator.

%
%

\subsection{Characteristic function estimator on $SO(3)$}

The characteristic function estimator presented here was used for example in \cite{Said2010} to provide a non-parametric estimation of the density of a compound Poisson process on the rotation group $SO(3)$. Recall that our notation is $\zeta^{\beta}_K (.)$ for this estimator. It is shown here that this estimator is in fact a kernel estimator and use can be made of the MISE expressions given in Section \ref{subsec:BVK}. Note that we consider the truncated version of the estimator, with $\ell_{max}=L$, given by:
\begin{equation}
\zeta^{\beta}_K (x) =  \sum_{\ell = 0}^{L}\sum_{n=-\ell}^{\ell}\sum_{m=-\ell}^{\ell} \left(2\ell+1\right) \widetilde{\hat{f}^{\ell}_{nm}} D^{\ell}_{nm}(x)
\end{equation}
 and where:
\begin{equation}
\widetilde{\hat{f}^{\ell}_{nm}}=\frac{1}{K}\sum_{k=1}^K \overline{D^{\ell}_{nm}}(X_k)
\end{equation}
The use of truncated version of the Fourier expansion is necessary for obvious computational reasons. Such a truncation has effects on the performance of the estimator as it only converges to a low-resolution version of the density. However, it is well-adapted to naturally band-limited functions. In the simulation in Section \ref{sec:simu}, we will investigate the influence of the bandwidth on the estimation performances. 

The definition of the characteristic kernel estimator leads by simple calculation to: 
\begin{equation}\label{eq:charfunctest}
\zeta^{\beta}_K (x) = \frac{1}{K} \sum_{k=1}^K \sum_{\ell \geq 0}(2\ell+1) \chi^{\ell}(X_k^{-1}x)
\end{equation}

This can be obtained thanks to the following property (see \cite{Hielscher2010} for example):
\begin{equation}
\sum_{n=-\ell}^{\ell} \sum_{m=-\ell}^{\ell} \overline{D^{\ell}_{nm}(x)} D^{\ell}_{nm}(y) = \chi^{\ell}(x^{-1}y)
\end{equation}
Note that this estimator indirectly estimates the density, as it is designed to estimate its Fourier coefficients. The estimator $\widetilde{\hat{f}^{\ell}_{nm}}$ is the characteristic function estimator. From equation (\ref{eq:charfunctest}), one can see that the characteristic function estimator is a Kernel estimator, with the specific kernel:

$$\Xi (.) = \sum_{\ell \geq 0}^{\ell_{max}=L} (2\ell+1) \chi^{\ell}(.)$$

From the expression of the MISE for a kernel estimator given in Equation (\ref{eq:MISE_K}), one can directly deduce the MISE for the characteristic function estimator. It suffices to replace the Fourier coefficients $\hat{\Xi}^{\ell}$, remembering the linearity property of the Fourier expansion, by the sum of the unit coefficients up to $\ell_{max}=L$. This make the characteristic function estimator a very ``simple'' kernel with constant Fourier coefficients. Illustration of its behaviour compared to other estimators will be presented in Section \ref{sec:simu}.

%
%

\subsection{Linear diffusive wavelet estimator}

Using the diffusive wavelet introduced in Section \ref{subsec:WTSO3}, it is possible to build an estimator of the density $f$ through its wavelet expansion. Recall that with our notation, the wavelet estimator is denoted $\zeta^{\gamma}_K (.,.)$. Note also that an extra scalar parameter is introduced when using the wavelet estimator. This is the scaling parameter of the wavelet transform. The estimator based on the wavelet transform consists in replacing the wavelet transform $WT_f(\rho,x)$ in the inversion formula given in Equation (\ref{eq:invWTSO3}) by its empirical estimate, obtained from the data sample. The wavelet coefficients are thus estimated using: 
\begin{equation}
\widetilde{WT}_f(\rho,x) = \frac{1}{K} \sum_{k=1}^K \check{\Psi}_{\rho}(X_k^{-1} x)
\end{equation}
The estimated density $\zeta^{\gamma}_K (x)$ takes the following expression when plugging the estimated coefficients in the inversion formula:
\begin{equation}
\zeta^{\gamma}_K (x,t) = \frac{1}{K}\sum_{k=1}^K \int_{t}^{+\infty}\int_{SO(3)} \check{\Psi}_{\rho}(X_k^{-1} y) \Psi_{\rho}(y^{-1} x) d\mu(y)\alpha(\rho)d\rho
\end{equation}
The 'scaling' coefficient $\rho$ is a parameter for this estimator. Ideal range should be $0$ to $+\infty$. Obviously it is not possible to reach the upper limit nad one can only numerically tend to very high values. As can be seen in the Simulation section, it is not a limitation for the use of the estimator. In the sequel, we will only keep the $t\geq 0$ constraint and keep in mind the upper limit issue.

It can easily be shown that the wavelet estimator can be expressed using the kernel $\kappa_t$ previsouly introduced in Equation (\ref{eq:heatkernelSO3}). 
The following equality thus stands:
\begin{equation}
\begin{array}{rcl}
\displaystyle{\int_{t}^{+\infty}\int_{SO(3)} \check{\Psi}_{\rho}(X_k^{-1}g) \Psi_{\rho}(g^{-1} x) d\mu(g)} & = & \displaystyle{ \int_{t}^{+\infty}\int_{SO(3)} \check{\Psi}_{\rho}(g') \Psi_{\rho}(g^{'-1} X_k^{-1}x) d\mu(g') }\\
& = & \displaystyle{ \int_{t}^{+\infty} \left( \check{\Psi}_{\rho} * \Psi_{\rho}\right) (X_k^{-1}x) \alpha(\rho)d\rho} \\
& = & \displaystyle{\kappa_t(X_k^{-1}x)}
\end{array}
\end{equation}
which is obtained by simple change of variable $g'=X_k^{-1}g$ and using the fact that the Haar measure is bi-invariant, inducing that $d\mu(X_k g')=d\mu(g')$. As a consequence, the linear wavelet estimator is:

\begin{equation}
\displaystyle{\zeta^{\gamma}_K (x,t) = \frac{1}{K}\sum_{k=1}^K \kappa_t(X_k^{-1}x)}
 \end{equation} 
where one can obviously see that this estimator is a kernel estimator. In order to study this estimator, one can look at its MISE with the bias and variance term. The bias term is: 

\begin{equation}
\begin{array}{rcl}
\displaystyle{ {\mathbb E}\left[\zeta^{\gamma}_K (x,t)\right]} & = & \displaystyle {{\mathbb E}\left[\frac{1}{K}\sum_{k=1}^K \int_{t}^{+\infty}\int_{SO(3)} \check{\Psi}_{\rho}(X_k^{-1} y) \Psi_{\rho}(y^{-1} x) d\mu(y)\alpha(\rho)d\rho\right]} \\
 & = & \displaystyle{ \frac{1}{K}\sum_{k=1}^K \int_{t}^{+\infty}\int_{SO(3)} \int_{SO(3)} f(z) \check{\Psi}_{\rho}(z^{-1} y) \Psi_{\rho}(y^{-1} x)d\mu(z) d\mu(y)\alpha(\rho)d\rho} \\
  & = & \displaystyle{ f * \int_{t}^{+\infty} \left( \check{\Psi}_{\rho}*\Psi_{\rho} \right) (x) \alpha(\rho) d\rho } \\
  & = & \displaystyle{\left(f*\kappa_t\right)(x)}
\end{array}
\end{equation}
Note that we have:
\begin{equation}
\lim_{t \rightarrow 0} {\mathbb E}\left[\zeta^{\gamma}_K (x,t)\right] = f(x)  
\end{equation}
Note that this is due to the fact that we used an approximate convolution identity. Thus, the bias of the linear wavelet estimator takes the following form: 
\begin{equation}
b^2_{\gamma,t}=\left\|f - {\mathbb E}\left[\zeta^{\gamma}_K(t)\right]\right\|_2^2 = \left\|f - f*\kappa_t \right\|_2^2 
\end{equation}
It is noticeable that, asymptotically with $t$, the bias term vanishes:
\begin{equation}
\lim_{t \rightarrow 0} b^2_{\gamma,t} = 0
\end{equation}
thanks to the definition of $\kappa_t$. This makes the wavelet estimator an unbiased estimator when $t$ reaches $0$. Now, for the variance term, notice first that:
\begin{equation}
\begin{array}{rcl}
\displaystyle{{\mathbb E}\left[ \zeta^{\gamma}_K(x,t) - {\mathbb E}\left[\zeta^{\gamma}_K(x,t)\right]\right]^2} & = & \displaystyle{var\left[\zeta^{\gamma}_K(x,t)\right]} \\
& = & \displaystyle{\frac{1}{K^2}  \sum_{k=1}^K \operatorname{var} \left[ \ \kappa_t(X_k^{-1}x) \right]} \\
 \end{array}
 \end{equation}


It comes then naturally that the variance term takes the form:
\begin{equation}
\begin{array}{rcl}
\displaystyle{\sigma^2_{\gamma,t}} & = & \displaystyle{\int_{SO(3)}var\left[\zeta^{\gamma}_K(x,t)\right] d\mu(x)}\\
 & = & \displaystyle{\frac{1}{K^2}\sum_{k=1}^K\int_{SO(3)} \left({\mathbb E}\left[\kappa_t^2(X_k^{-1}x)\right] - {\mathbb E}^2\left[\kappa_t(X_k^{-1}x)\right] \right) d\mu(x)} \\
  & = & \displaystyle{ \frac{1}{K^2}\sum_{k=1}^K\int_{SO(3)} \int_{SO(3)} \kappa_t^2(z^{-1}
  x)f(z)d\mu(z)d\mu(x) - \frac{1}{K^2}\sum_{k=1}^K \int_{SO(3)} \left(f*\kappa_t\right)^2(X_k)d\mu(X_k)}
 \end{array}
\end{equation}
where we made use of the bi-invariance of the Haar measure and of the independence between the samples $X_k$. Finally, it is possible to give the expression of the variance as:
\begin{equation}
\sigma^2_{\gamma,t}={\mathbb E}\left[ \left\| \zeta^{\gamma}_K (x,t) - {\mathbb E}\left[\zeta^{\gamma}_K (x,t)\right] \right\|_2^2 \right] = \frac{1}{K} \left( \left\| \kappa_t \right\|_2^2 - \left\| f*\kappa_t \right\|_2^2\right)
\end{equation}





As was noticed earlier, the wavelet estimator is a kernel estimator. Using the general formula for the MISE of kernel estimators and the expressions of the bias and variance given above, one gets the following expression for the MISE:


\begin{equation}
MISE(\zeta^{\gamma}_K,t)=\sum_{\ell \geq 1} \left( (2\ell +1) \left.\hat{f}^{l}\right.^2  (1 - e^{-\ell(\ell+1)t})^2  + \frac{(2\ell +1)}{K} e^{-2\ell(\ell+1)t} \left(1 - \left.\hat{f}^{l}\right.^2\right)\right)
\end{equation}

It is then possible to give a first order approsimation of this expression. This takes the form:
\begin{equation}
MISE(\zeta^{\gamma}_K,t)=\sum_{\ell \geq 1} \left[\frac{(2\ell +1)(1-2\ell(\ell+1)t)}{K} +  \left.\hat{f}^{l}\right.^2  \frac{(2\ell+1)}{K} \left( 2\ell\left(\ell+1\right)t\left(1+K\right)-1\right) \right]
\end{equation}
This expression will be used in the Simulation Section. It is of interest to note that the behaviour of the MISE with parameter $t$ is as follows:
\begin{equation}
\lim_{t \rightarrow 0} MISE \left(\zeta^{\gamma}_K,t \right) = \sum_{\ell \geq 1}\frac{(2\ell + 1 )}{K}(1-\left.\hat{f}^{l}\right.^2)
\end{equation}

Also notice that as $f$ is a probability density, then $f(x) \geq 0$ almost everywhere and we have that $\hat{f}^{0}=1$ because $\hat{f}^{0}=\int_{SO(3)}f(x)d\mu(x)$.

We have presented the MISE for the kernel estimator and for the characteristic function and linear wavelet estimators. The three studied estimators belong to the kernel estimator family, but they do have differences. The given expressions will allow allow their behaviours to be compared. The following section highlights these differences.  




\section{Simulations}
\label{sec:simu}
In this section, we present some simulations that illustrate the differences between the three considered estimators studied in this paper. The comparison is performed in terms of the respective MISE computed in the context of estimation of a mixture of densities on the rotation group $SO(3)$.

\subsection{Definition of a test function}
 First, we have to define a test function $f$ on the rotation group, which we take similar to the one studied in \cite{Hielscher2013}. This mixture takes the form:
\begin{equation}
f(\bm x) = 0.2 + 0.7\psi_{\text{VP}}^{30}(\bm R(\bm e_1,\pi/6)\cdot \bm x) + 0.1\psi_{\text{VP}}^{45}(\bm R(\bm e_2,4\pi/9)\cdot \bm x),
\end{equation}
where $\bm R(\eta, \omega)$ denotes the rotation of angle $\omega$ and axis $\eta$, and $\psi_{\text{VP}}^\kappa$ is the de La Vallee Poussin kernel. This kernel has a closed form expression, which reads
\begin{equation}
\psi_{\text{VP}}^\kappa(x)= \frac{(2\kappa+1)2^{2\kappa}}{\binom{2\kappa+1}{\kappa}}\cos^2\left(\frac{\omega(x)}{2}\right) = \binom{2\kappa+1}{\kappa}^{-1}\sum_{\ell=0}^\kappa(2\ell+1)\binom{2\kappa+1}{\kappa-\ell}\chi^\ell(x).
\end{equation}
The parameter $\kappa$ can be understood as the analog of the $\rho$ parameter in the heat kernel wavelet approah, that is it plays the equivalent role of a bandwidth. 
\subsection{Computing the MISE}
We recall that the MISE can be expressed in the Fourier domain for a kernel function $\Xi$ by
\begin{equation}
MISE(\Xi) = \sum_{\ell \geq 1} \left( (2\ell +1) \left.\hat{f}^{\ell}\right.^2  (1 - \left.\hat{\Xi}^{\ell}\right.)^2  + \frac{(2\ell +1)}{K} \left.\hat{\Xi}^{\ell}\right.^2 \left(1 - \left.\hat{f}^{\ell}\right.^2\right)\right).
\end{equation}
The last expression of the MISE  in terms of the Fourier coefficients of both the test function $f$ and the chosen kernel $\Xi$ allows us to compute the MISE in a simple way. The Fourier coefficients of the kernel are in general known, as it is the case here for the kernel used (de la Vallée Poussin, Heat kernel, and characteristic function kernel). However the coefficients of $f$ have to be computed with a numerical implementation of the Fourier transform on the rotation group.

We used in our simulation an implementation of the FFT on $SO(3)$ as proposed by Kostelec and Rockmore \cite{kostelec2008ffts}. This FFT is based on a equiangular sampling of the Euler angles. Since our test function $f$ is a linear combination of de la Vall\'ee Poussin kernels, it is bandlimited by the largest $\kappa$ value chosen, that is in our case $\kappa = 45$. The FFT was thus performed up to degree $L = 49$, leading to 166650 complex valued Fourier coefficients $\hat{f}_{lmn}$. The energy per degree $\hat{f}^{\ell^2}$ can be computed like:
\begin{equation}
\hat{f}^{\ell^2} = \frac{1}{(2\ell+1)}\sum_{m,n = -\ell}^\ell\vert f_{lmn}\vert^2. 
\end{equation}
This expression, up to the maximum value of $\ell$ is then plugged into the MISE expression.
\subsection{Results}
\begin{figure}
\centering
\includegraphics[scale=.5]{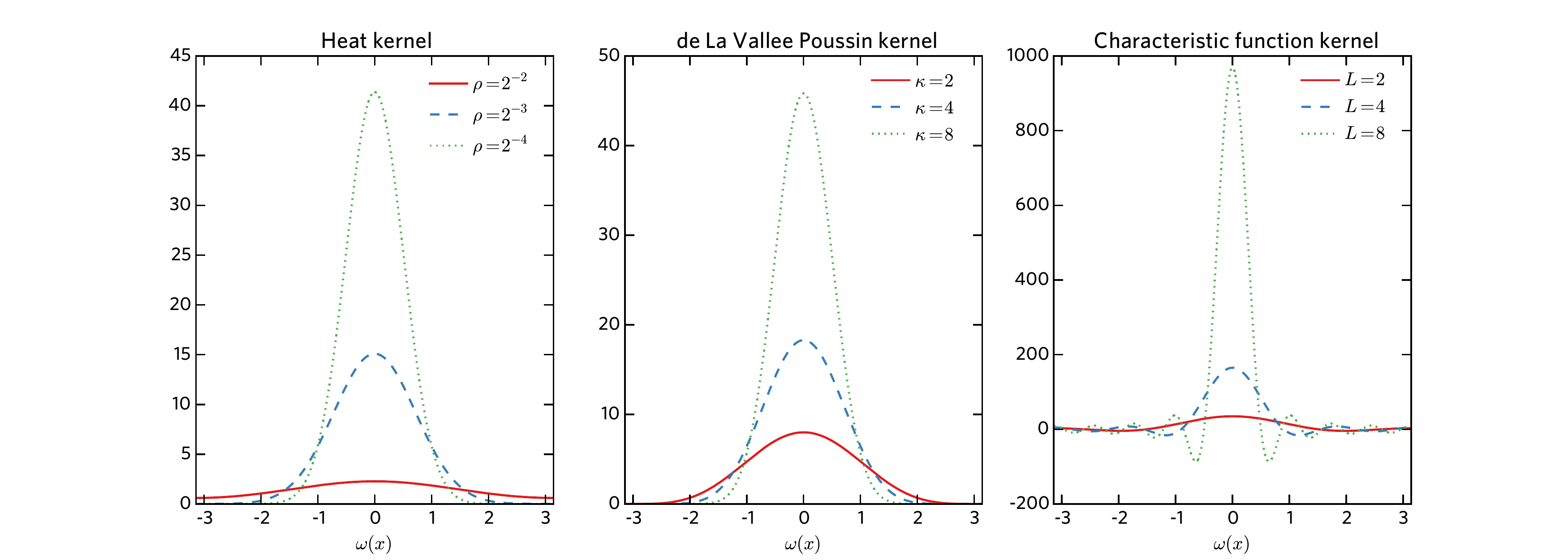}\caption{Kernels used in the simulations. From left to right: Heat kernel, La Vall\'ee Poussin kernel and characteristic function kernel. Kernels are displayed for different values of the ``bandwidth'' parameter.}\label{fig:kernels}
\end{figure}

In figure \ref{fig:kernels} we have plotted the kernels used in our simulations. We consider here only three kernel types, the de La Vall\'ee Poussin kernel, the Heat kernel (sometimes called the Gauss-Weierstrass kernel) and the characteristic function kernel. For each of these kernels we looked at different bandwidth parameters, respectively $\kappa$, $\rho$, and $L$, the latter being the truncation order in the characteristic function estimator. The plots in figure \ref{fig:kernels} emphasize the role of the respective bandwidth, as the concentration of the kernel functions increases with larger bandwidths. We also note that the de La Vall\'ee Poussin and heat kernels are nonnegative kernels, whereas the characteristic function kernel exhibits negative values. Such behaviour has drawbacks, especially when the density to estimate exhibits narrow modes.

\begin{figure}
\centering
\includegraphics[scale=.5]{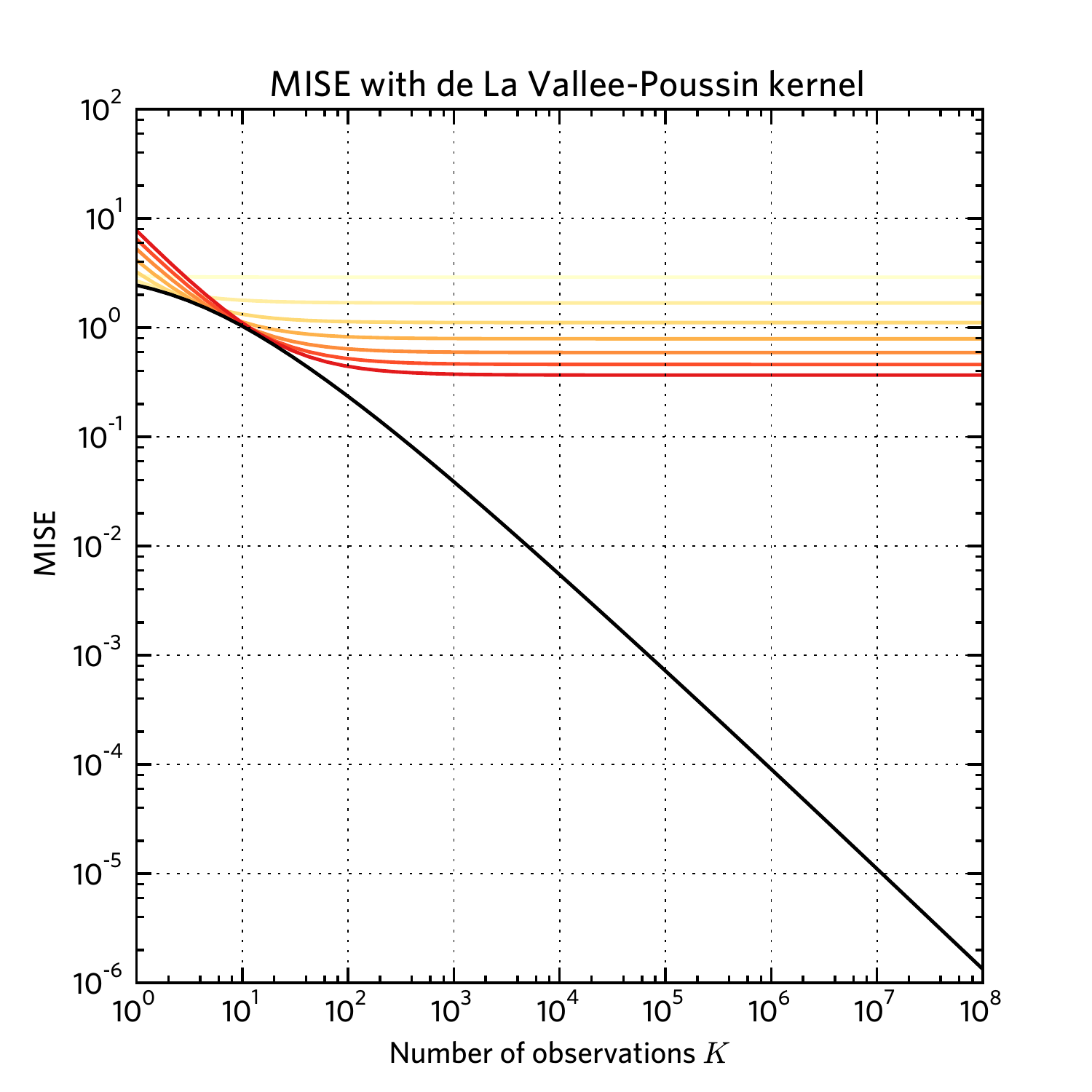}\caption{MISE computed for a de la Vallee Poussin kernel type, for bandwidth values $\kappa  = 1,8,17,22,29,36;43$ (from light yellow to dark red). The optimal MISE bound given in Equation (\ref{eq:MISE_optimal_bound}) is displayed in black for comparison. }\label{fig:MISE1}
\end{figure}

\begin{figure}
\centering
\includegraphics[scale=.5]{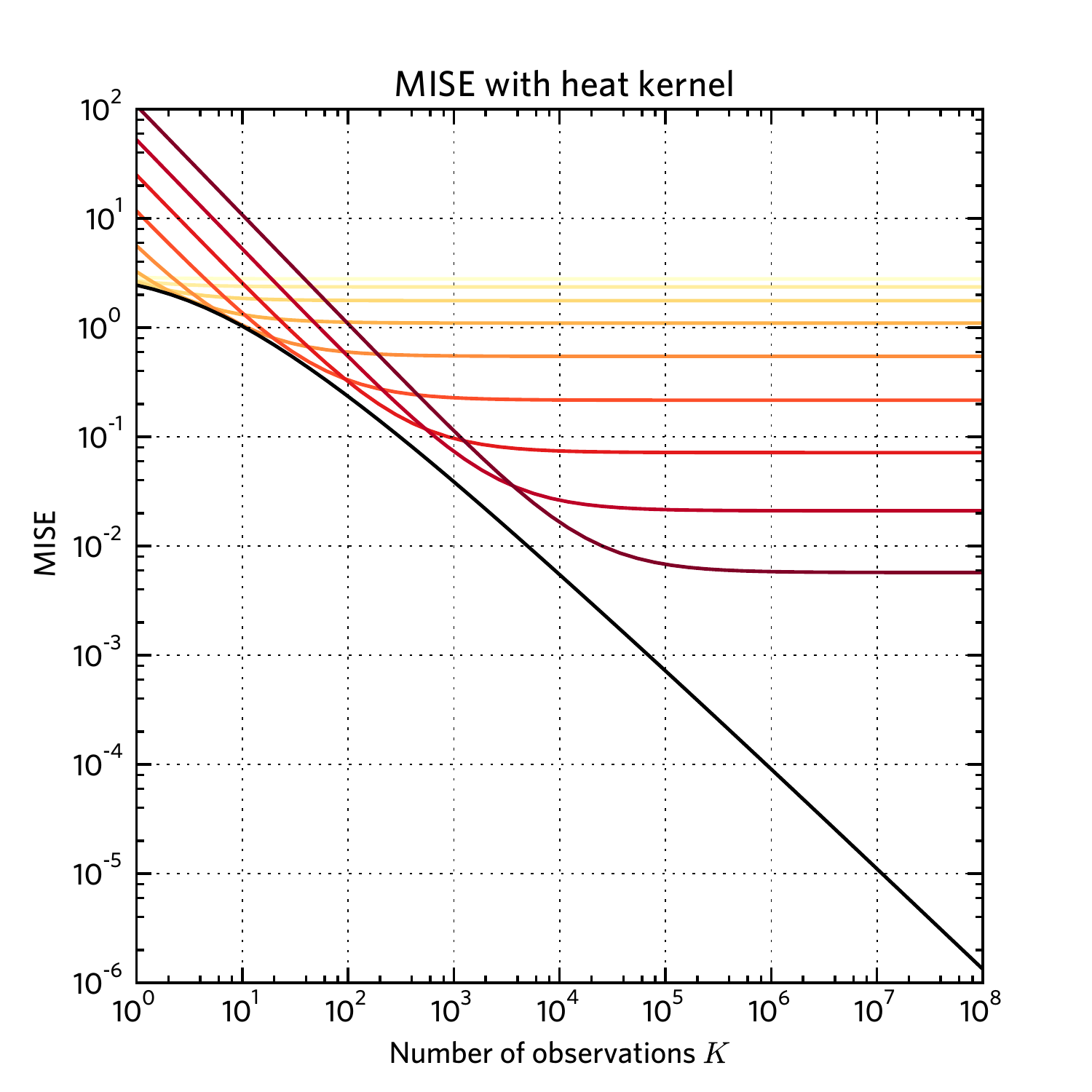}\caption{MISE computed for the heat kernel, for bandwidth values $\rho  = 2^{-j}$ with $j=1,2,\ldots,9$ (from light yellow to dark red). The optimal MISE bound given in Equation (\ref{eq:MISE_optimal_bound}) is displayed in black for comparison. }\label{fig:MISE2}
\end{figure}

\begin{figure}
\centering
\includegraphics[scale=.5]{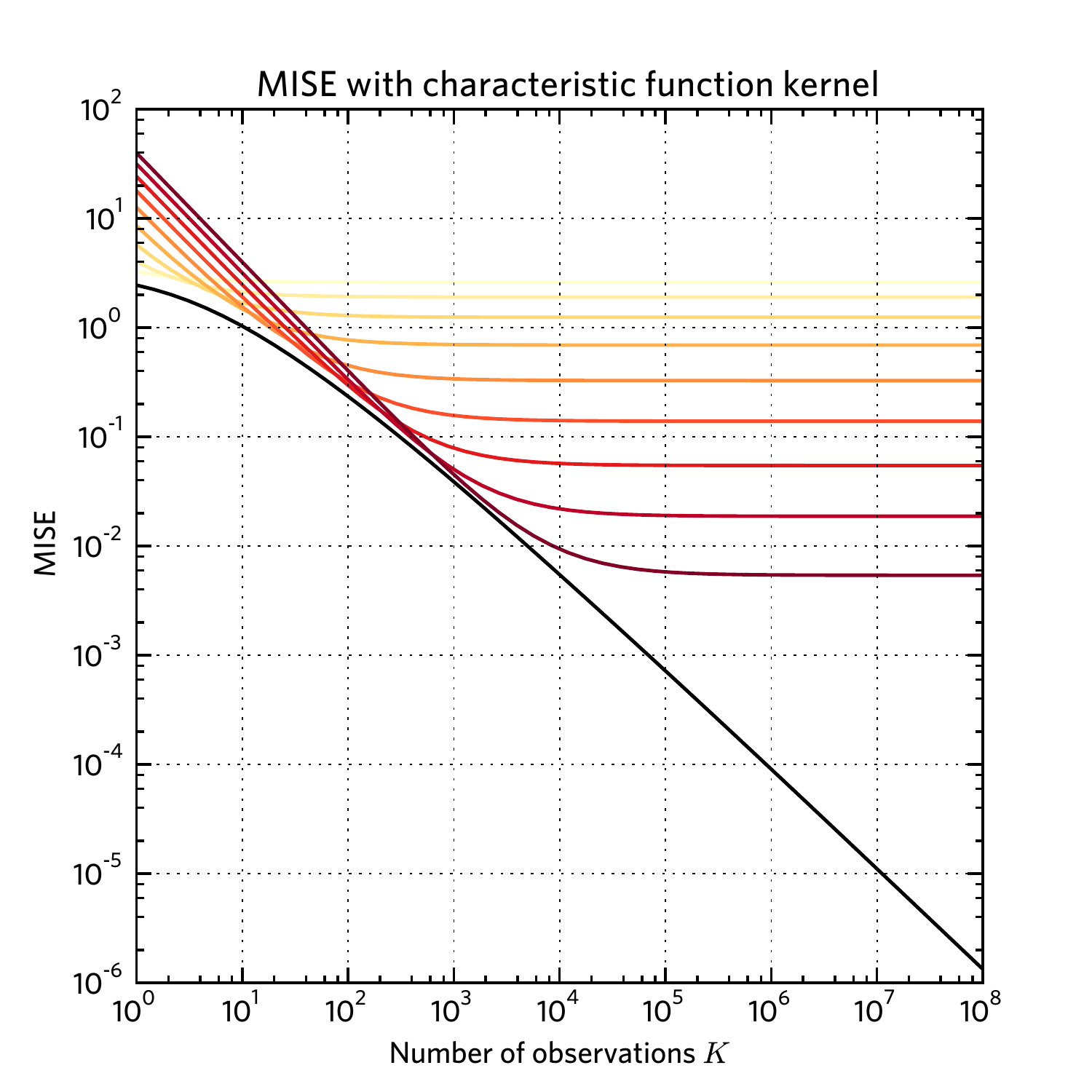}\caption{MISE computed for the characteristic function kernel, for bandwidth values $L=1,2,\ldots,9$ (from light yellow to dark red).The optimal MISE bound given in Equation (\ref{eq:MISE_optimal_bound}) is displayed in black for comparison.	}\label{fig:MISE3}
\end{figure}

Before stepping into the MISE figures, we recall an important result from Hielscher \cite[theorem 3]{Hielscher2013}, which gives the MISE optimal bound for the function $f$. It reads explicitly as 
\begin{equation}\label{eq:MISE_optimal_bound}
MISE_\text{opt} = \sum_{\ell=1}^\infty (2\ell+1)\frac{\hat{f}^{\ell^2}(1-\hat{f}^{\ell^2}}{(K-1)\hat{f}^{\ell^2}+1}
\end{equation}
where we have adapted the prefactor in $(2\ell +1)$ due to our normalization choice. This optimal MISE is displayed on figures \ref{fig:MISE1}, \ref{fig:MISE2} and \ref{fig:MISE3} for comparison with the studied estimators.

We investigate now the behavior of each kernel and its influence on the MISE for different values of bandwidth parameters. It is important to note that in each experiment, bandwidths were fixed, as opposed to what was done in \cite{Hielscher2013} were the bandwidths were chosen as optimal with respect to the test function and the number of observations. Our approach is indeed motivated by a multiresolution approach, with the minimum amount of information about the function $f$ being incorporated in the kernels. 

In figure \ref{fig:MISE1}  we have the MISE evaluated for de la Vall\'ee Poussin kernel with bandwidth parameter ranging for $\kappa  = 1,8,17,22,29,36,43$ along the theoretical lower bound. In figure \ref{fig:MISE2}  we have the MISE evaluated for the heat kernel with bandwidth parameter ranging from $\rho = 2^{0}$ to $\rho = 2^{-9}$, (which corresponds to a dyadic scaling) along the theoretical lower bound. Finally, in figure \ref{fig:MISE3} we have the MISE evaluated for the characteristic function kernel with bandwidth parameter ranging from $L=1$ to $L=9$, along the theoretical lower bound.

As expected with a multiresolution approach, for a given bandwidth the MISE reaches a lower threshold after a sufficient number of observations $K$. Moreover as the bandwidth increases the threshold is lowered. Comparing the three kernel types, the de La Vall\'ee Poussin kernel seems to converge slower than the heat or characteristic function kernel, as the threshold lowers slowly. Therefore, the de La Vall\'ee Poussin kernel seems less appropriate for a multiresolution analysis.
The heat kernel and characteristic function kernel are performing better, with a slight advantage to the characteristic function kernel as it comes closest to the MISE optimal bound. We note however that the characteristic function kernel does not have the nonnegative property of the heat kernel, which may be critical when the density to estimate possesses sharpened modes. It is also noticeable that for small sizes of samples (from $10^1$ to $10^2$) the wavelet estimator gets closer to the optimal bound. This suggests that an adaptive choice of the parameter $\rho$ could lead to very good estimates in most of the situations.    



\section{Conclusion}

We have demonstrated that the characteristic function estimator and the linear Heat wavelet estimator on $SO(3)$ both belong to the larger family of kernel estimators for densities on $SO(3)$. The characteristic function estimator consists in a kernel estimator with constant Fourier coefficients up to a maximum bandwidth, while the wavelet estimator leads to a Heat kernel with coefficients driven by the scaling parameter of the heat wavelet family. The MISEs of the introduced estimator have been presented and illustration of the differences between heat kernel, characteristic function and De La Vall\'ee Poussin kernels investigated. The diffusive wavelet based estimator combines the nice property of converging faster than the De La Vall\'ee Poussin kernel and of being strictly positive (as opposed to the characteristic function kernel), allowing good performances in many configurations. The results presented demonstrate the advantage of using the heat wavelet transform for probability density estimation on the rotation group $SO(3)$. Future work should consist in studying the nonlinear (thresholded) version of the heat wavelet estimator which is known to perform better than linear wavelet estimator on the real line. Validation on real datasets should also be investigated to completely validate the proposed estimator. 




\bibliographystyle{IEEEtran}
\bibliography{IEEEabrv,DiffWav_SO3_biblio}







\end{document}